\documentclass[12pt]{amsart}
\usepackage[colorlinks=true, pdfstartview=FitV, linkcolor=blue, citecolor=blue, urlcolor=blue]{hyperref}
\usepackage{amssymb}
\calclayout

\def\C{\mathbb {C}}
\def\R{\mathbb {R}}

\def\N{\mathbb {N}}
\def\Z{\mathbb {Z}}

\def\inv{^{-1}}

\def\GL{\operatorname{GL}}

\def\phi{\varphi}

\def\codim{\operatorname{codim}}

\numberwithin{equation}{subsection}

\newtheorem{theorem}[subsection]{Theorem}
\newtheorem{lemma}[subsection]{Lemma}

\newtheorem{corollary}[subsection]{Corollary}

\theoremstyle{definition}

\theoremstyle{remark}

\newtheorem{remark}[subsection]{Remark}

\newtheorem{example}[subsection]{Example}

\title[Isomorphisms preserving invariants]{\boldmath Isomorphisms preserving invariants} 
 \author{Gerald W. Schwarz}
\thanks{Partially supported by NSA Grant H98230-06-1-0023}
\address{Department of Mathematics\\
Brandeis University\\
Waltham, MA 02454-9110}
\email{schwarz@brandeis.edu}

\begin{document}
 \begin{abstract}
 Let $V$ and $W$ be finite dimensional real vector spaces and let $G\subset\GL(V)$ and $H\subset\GL(W)$ be finite subgroups.
 Assume for simplicity that the  actions contain no reflections. Let $Y$ and $Z$ denote the real algebraic varieties corresponding to $\R[V]^G$ and $\R[W]^H$, respectively. If $V$ and $W$ are quasi-isomorphic, i.e., if there is a linear isomorphism $L\colon V\to W$ such that $L$   sends $G$-orbits to $H$-orbits and $L\inv$ sends $H$-orbits to $G$-orbits, then $L$ induces an isomorphism of $Y$ and $Z$. Conversely, suppose that $f\colon Y\to Z$ is a germ of a diffeomorphism sending the origin of $Y$ to the origin of $Z$. Then we show that $V$ and $W$ are quasi-isomorphic, This result is  closely related to a theorem of Strub \cite{Strub}, for which we give a new proof.  We also give a new proof of a result of \cite{KrieglLosikMichor03} on lifting of biholomorphisms of quotient spaces. 
\end{abstract}
\maketitle

\section{Introduction}
First some motivation. Let $G$ be a finite group and let $M$ and $N$ be smooth $G$-manifolds. We give the orbit spaces $M/G$ and $N/G$ a smooth structure by declaring that $C^\infty(M/G)=C^\infty(M)^G$ and $C^\infty(N/G)=C^\infty(N)^G$. We also have a stratification of $M/G$ by isotropy type, where the type of a point $Gm$ of $M/G$ is the conjugacy class of the isotropy group $G_m$ of $m\in M$.   Suppose that $M$ and $N$ are equivariantly diffeomorphic. Then  $M/G$ and $N/G$ are ``diffeomorphic''  and the diffeomorphism preserves isotropy type. Conversely, suppose that $M/G$ and $N/G$ are diffeomorphic by a strata preserving diffeomorphism $f$. Can we lift $f$ to an equivariant diffeomorphism of $M$ and $N$? First of all, $f$ must give a diffeomorphism of $M^G$ and $N^G$ (assuming these sets are nonempty). Let $v\in M^G$ and set $V=T_vM$. Let $f(v)=w\in N^G$ and set $W=T_wN$. Then the differentiable slice theorem says that there is a neighborhood of $v$ in $M/G$ diffeomorphic to $U\times V/G$ where $U$ is a neighborhood of $v$ in $M^G$. It follows that   $f$ induces a strata preserving diffeomorphism of $V/G$ with $W/G$. The first question to answer is if the existence of a diffeomorphism of $V/G$ with $W/G$ implies that $V$ and $W$ are isomorphic representations of $G$. The answer is ``yes,'' up to an automorphism of $G$, by a theorem of Strub \cite{Strub}. We generalize the problem by replacing the requirement that strata are preserved by a weaker condition, we allow for a different group $H$ to act on $W$   and we allow the diffeomorphism $f$ to be defined on the Zariski closures of the orbit spaces, in a sense made precise below. This ends our motivation, and the precise description of our results follows.

 Let $V$ and $W$ be finite dimensional real vector spaces and let $G\subset\GL(V)$ and $H\subset\GL(W)$ be finite subgroups. Let $Y$ and $Z$ denote the real algebraic varieties corresponding to $\R[V]^G$ and $\R[W]^H$, respectively.  We say that a linear isomorphism $L\colon  V \to W$ is a  \emph{quasi-isomorphism\/}  if $L$  sends $G$-orbits to $H$-orbits and  $L\inv$ sends $H$-orbits to $G$-orbits. We    say that $V$ and $W$ are \emph{quasi-isomorphic} if such an $L$ exists. If $L$ is a quasi-isomorphism, then $g\mapsto L\circ g\circ L\inv$ gives an isomorphism of $G$ and $H$ such that, with $G$ acting on $W$ via the isomorphism,   $L$ is equivariant.  

If $V$ and $W$ are quasi-isomorphic, then $Y$ and $Z$ are isomorphic. Conversely, we show that if $f\colon Y\to Z$ is a germ of a diffeomorphism sending the origin of $Y$ to the origin of $Z$ such that $f$ maps the closures of the  codimension one strata of $Y$   onto those of $Z$, then $V$ and $W$ are quasi-isomorphic. The origin of $Y$ is the image  of $0\in V$ in $Y$   and similarly for $Z$.

There are closely related results in the literature. Let $Y_0$ and $Z_0$ denote the images of $V$ and $W$ in $Y$ and $Z$, respectively. Then $Y_0$ is isomorphic to the orbit space $V/G$ and   $Z_0$ is isomorphic to $W/H$. Strub \cite{Strub} has shown that if   there is a diffeomorphism of neighborhoods of the origins in   $Y_0$ and $Z_0$, then   the $G$ and $H$ actions are quasi-isomorphic. Prill \cite{Prill} showed the analogue of Strub's result  in the case of complex representations $U$ of finite complex groups $K$ not containing pseudoreflections. Here one considers the complex variety $U/K$ corresponding to $\C[U]^K$. Gottschling \cite{Gottschling} and Prill \cite{Prill} proved that biholomorphisms of $U/K$ lift to equivariant biholomorphisms of $U$, again assuming that $K$ contains no pseudoreflections.  
  In \cite{Losik01} and \cite{KrieglLosikMichor03} there are results about lifting of isomorphisms of   complex quotients $U /G\simeq U' /H$ where $U$ and $U'$ are complex $G$ and $H$-modules, respectively. Lyashko \cite{Lyashko} showed that lifting holds in case $G=H$ is a Weyl group and $U=U'$    is the complexification of the standard real representation. Losik \cite{Losik01} proves  that diffeomorphisms of $Y_0$ lift to diffeomorphisms  of $V$ which preserve $G$-orbits. 
  We give  new proofs  of the lifting result of \cite{KrieglLosikMichor03} and  of  the main result of \cite{Strub}

We thank K. Abe for posing  questions leading to this note and for making us aware of \cite{Strub}.

\section{Isomorphisms of spaces of invariants}

 Let   $p_1,\dots,p_m$  be homogeneous generators of $\R[V]^G$  and   let $q_1,\dots,q_n$ be homogeneous generators  of $\R[W]^H$. Let $d_j$  denote the degree of $p_j$, $j=1,\dots,m$ and let $e_j$ denote the degree of $q_j$, $j=1,\dots,n$. Let $p=(p_1,\dots,p_m)\colon V\to \R^m$ and $q=(q_1,\dots,q_n)\colon W\to \R^n$. Let $Y_0$ denote $p(V)$ and let $Z_0$ denote $q(W)$. Then $Y_0$ and $Z_0$ are closed and semialgebraic where $Y_0\simeq V/G$ and $Z_0\simeq W/H$. We may identify $Y$ and $Z$ with the Zariski closures of $Y_0$ in $\R^m$ and of $Z_0$ in $\R^n$.
Since the $p_i$ are homogeneous,  $Y_0$ is stable under the $\R^*$ action on $\R^m$ which sends $t\in\R^*$ and $y= (y_1,\dots,y_m)\in \R^m$ to $t\cdot y:=(t^{d_1}y_1,\dots, t^{d_m}y_m)$. It follows that $Y$ is stable under this action. Similarly we have an action of $\R^*$ on $\R^n$ which preserves $Z_0$ and $Z$.  
Let $V_\C=V\otimes_\R\C$ denote the complexification of $V$. Then $p\colon V_\C\to\C^m$ has image the complex variety $Y_\C$ corresponding to $\C[V_\C]^G$. Moreover,  $Y$ is $Y_\C\cap\R^m$. The  strata of $Y$ and $Y_0$ are the collections of points whose preimages in $V_\C$ have conjugate isotropy groups. There are clearly finitely many strata. We say that a stratum is \emph{codimension one\/}  if it is of codimension one in $Y_0$ or $Y$, as the case may be. Equivalently, the isotropy groups corresponding to the stratum have order 2 and are generated by   reflections.   The \emph{principal points\/} of $Y$ and $Y_\C$ are those points whose inverse images in  $V_\C$ have trivial isotropy groups.

\begin{example}\label{ex:R1}
Let $V=\R$ and $G=\{\pm 1\}$ acting by multiplication. Then $p=x^2\colon \R\to\R$. We have that $Y_0=\R^+$ and $Y=\R$. The codimension one stratum  of $Y$  is the origin  and the principal points are $Y\setminus\{0\}$.  
\end{example}

\begin{example}\label{ex:R2}
Let $V=\R^2$ and $G=\Z/k\Z$ acting via rotations, $k\geq 2$. Let $x$ and $y$ be the usual coordinate functions on $V$ and set $z=x+iy$. Then the $G$-invariant polynomials are generated by $z\bar z$ and the real and imaginary parts of $z^k$. Thus $p_1=x^2+y^2$, $p_2=x^k-\binom k2 x^{k-2}y^2+\dots$ and $p_3=kx^{k-1}y-\binom k3 x^{k-3}y^3+\dots.$  generate $\R[V]^G$. They  satisfy the relation $p_2^2+p_3^2=p_1^k$. Hence $Y=\{(y_1,y_2,y_3)\in \R^3\mid y_1^k=y_2^2+y_3^2\}$ and one can show that $Y_0=\{(y_1,y_2,y_3)\in Y\mid y_1\geq 0\}$. If $k$ is odd, then $Y=Y_0$, while $Y$ and $Y_0$ differ if $k$ is even.
Here the strata are again $\{0\}$ and $Y\setminus\{0\}$, so there are no codimension one strata.
\end{example}
 
 Assume that we have a germ of a diffeomorphism $f\colon Y\to Z$  sending $0$ to $0$. This means that $f$ extends to a   smooth germ of a mapping from $\R^m$ to $\R^n$ sending $Y$ to $Z$ and similarly for $f\inv$. We assume that  $f$ sends the closures of the codimension one strata of $Y$ onto the closures of the codimension one strata of $Z$.      
 \begin{lemma} \label{lem:principal}
   Let $f\colon  Y\to Z$ be a local diffeomorphism as above. Then $f$ maps  principal points of $Y$ to principal points of $Z$.
   \end{lemma}
   
   \begin{proof} Everything we do should be understood to happen in neighborhoods of $0$ in $Y$ and $Z$ where $f$ and $f\inv$ are defined. Let $C$ and $D$ be the closures of the codimension one strata in $Y$ and $Z$, respectively. It follows from the Shepherd-Todd-Chevalley theorem that $y\in Y$ is a smooth point if and only if the isotropy group of a preimage $y_0\in V_\C$ is  trivial or is  generated by reflections. In the latter case let  $r$ be one of the reflections in $G_{y_0}$. Then $y$ lies in the closure of the codimension one stratum consisting of points whose isotropy group is generated by $r$. Thus a point in   $Y\setminus C$ is smooth in $Y$  if and only if it is principal. Hence the principal stratum of $Y$ consists of the smooth points not in $C$, so that $f$ induces an isomorphism of the principal points of $Y$ and $Z$.   
   \end{proof}

   Let $y\in Y$ be close to zero so that $f$ is defined on $(0,1]\cdot y$. For $t\in (0,1]$, let $f_t$ denote $t\inv\cdot f(t\cdot y)$.  We say that $f$ is \emph{quasilinear\/} if $f_t=f$ for  all $t\in (0,1]$. In this case, of course, $f$ is the restriction of a global diffeomorphism of $Y$ and $Z$ which is induced by a polynomial mapping.

 \begin{lemma} \label{lem:f_0} Let $f\colon Y\to Z$ be a local diffeomorphism as above.   Then $f_t(y)$ converges uniformly to a limit $f_0(y)$  as $t\to 0$ for $y$ in a neighborhood of $0$ in $Y$. The mapping $f_0$ is   a quasilinear polynomial isomorphism of  $Y$ and $Z$ which preserves the closures of the codimension one strata.
   \end{lemma}
   
    \begin{proof}    There are invariant inner products on $V$ and $W$  and we may assume that $p_1$ and $q_1$ are the corresponding quadratic forms.  Let $Y'$ and $Z'$ denote the principal points of $Y$ and $Z$, respectively. Then $Y'$ and $Z'$ are $\R^*$-stable and have finitely many components. The image $Z_1\subset Z$ of the principal points of $W$ is connected and open in $Z$ and it is $\R^*$-stable.  Since $q(W)=Z_0$ is closed in $Z$, any limit point of $Z_1$ which is not in $Z_1$ is not a principal point. Thus $Z_1$ is a connected component of $Z'$.  It follows that  $f\inv Z_1$ lies in and  contains a neighborhood of $0$ of a connected component $Y_1$ of $Y'$.   
 Let $y\in Y_1$ such that $f$ is defined on $(0,1]\cdot y$. Set $c(t)=f(t\cdot y)$. Then $c(t)=(c_1(t),\dots,c_n(t))$ is a curve in $Z_0$ and from  \cite[Lemma 2.1]{Losik01} we see that $c_i(t)$ vanishes at least to  order $e_i$ at $t=0$. Thus $f_t(y)$ converges  to a limit $f_0(y)$. Now consider the Taylor polynomial $Tf$ of $f$ up to the maximum degree $s$ of the $q_i$. Then the $i$th component $T_i(f)$ is a sum of monomials $\sum c_\alpha y^\alpha$ where $c_\alpha\in \R$ and $\alpha=(\alpha_1,\dots,\alpha_m)\in\N^m$. Let $|\alpha|$ denote $d_1\alpha_1+\dots+d_m\alpha_m$. Then our calculation shows that the sum $\sum_{|\alpha|< e_i}c_\alpha y^\alpha$ vanishes on $Y_1$. Since $Y_1$ is   Zariski dense in  $Y$    the sum vanishes on $Y$.  We may change $f$ without changing its restriction   to $Y$ so  that $T_i(f)$ has nonzero coefficients only for monomials $y^\alpha$ with $|\alpha|\geq e_i$. Then the Taylor polynomial of $f$ up to degree $s$ induces the mapping $f_0$ and the convergence of $f_t$ to $f_0$ is uniform for $y\in Y$ in a neighborhood of $0$.     Since $f_0$ is a limit of  maps preserving the closures of the codimension one strata, $f_0$ does also, and $f_0$ is clearly   quasilinear. It is also invertible, with inverse $(f\inv)_0$. 
   \end{proof}

   \begin{remark}
   It is clear that $f_0$ induces   a complex   isomorphism of $Y_\C$ and $Z_\C$, also denoted $f_0$, preserving the closures of the codimension one strata. This follows from the fact that the Zariski closures of the strata of $Y_0$ in $Y_\C$ are the strata of $Y_\C$, and this correspondence preserves codimension \cite[5.8]{SchLifting}.
      \end{remark}
 
 \begin{remark}
 In Examples  \ref{ex:R1}  and \ref{ex:R2}  (with $W=V$ and $H=G$) the mapping $f_0$ is linear.  This is obvious in Example \ref{ex:R1}  since the generators are all quadratic. Example \ref{ex:R2} requires more work.
  \end{remark}

   \begin{corollary} Let $f_0\colon Y_\C\to Z_\C$ be a quasilinear isomorphism which induces an isomorphism of the closures of the codimension one strata. Then there is a quasi-isomorphism $F\colon V_\C\to W_\C$ which induces  $f_0$.  
  \end{corollary}
  
   \begin{proof}
  The proof of Lemma \ref{lem:principal} shows that $f_0\colon Y_\C\to Z_\C$ induces an isomorphism of  the principal points. Since it also preserves the closures of the codimension one strata, it has a biholomorphic   lift $F\colon V_\C\to W_\C$ such that $F$ and $F\inv$ preserve orbits  \cite{KrieglLosikMichor03}.  
  Let $F_t(v)=t\inv F(tv)$ for $t\in\R^*$ and $v\in V$. Then $F_t$ is a continuous family of biholomorphic lifts  of $f_0$, so that we must have $F_t=F$ for all $t$ and hence $F=\lim_{t\to 0} F_t$ is linear. Thus $F$ is a quasi-isomorphism.
  \end{proof}
   
   Let $h\in H$. Set $W_h:=W_{h+}\oplus iW_{h-}$ where $W_{h\pm}$ denotes the $\pm 1$-eigenspace of $h$.
   
\begin{lemma} \label{lem:wg}
Let $z'\in Z$. Then there is an $h\in H$ such that   $z'\in p(W_h)$. 
\end{lemma}

\begin{proof}
Let $z\in W_\C$ such that $q(z)=z'$. Since the polynomials making up $q$ have real coefficients, $q(z)=\overline{q(z)}=q(\bar z)$. Since $q(z)=q(\bar z)$, there is an $h\in H$ such that $\bar z = hz$. Then $z=w_1+iw_{-1}$ for some $w_1\in W_{h+}$ and $w_{-1}\in W_{h-}$.   
\end{proof}

\begin{theorem}\label{thm:main}
Let $f\colon Y\to Z$ be a local diffeomorphism  sending $0$ to $0$ which preserves the closures of the codimension one strata. Then there is a quasi-isomorphism $F\colon V\to W$.
\end{theorem}

\begin{proof}
We have shown that there is a  quasi-isomorphism $F\colon V_\C\to W_\C$ which induces $f_0$ on $Y$. Thus $F(V)$ has to be an $H$-stable linear subspace of $W_\C$ which maps to $Z$. Since $H$ is finite, it follows from Lemma \ref{lem:wg} that $F(V)=W_h$ for some $h\in H$ such that $h^2=1$ and such that $h$ is central in $H$ (else $F(V)$ is not $H$-stable). Since $W_h$ is $H$-equivariantly isomorphic to $W$, we see that   $V$ and $W$ are quasi-isomorphic.
\end{proof}

\begin{remark}
Let  $F$, $h$ and $f$ be as in the proof above. The linear mapping from $\R[W]^H$ to $\R[W]^H$ which sends $f(w_1+w_{-1})$ to $f(w_1+iw_{-1})$ induces a graded automorphism of $\R[W]^H$, hence a ``linear'' automorphism $\phi$ of $Z$.  Then $\phi$ identifies  $Z_0$ and $q(W_h)$. Our  diffeomorphism $f$  is just   a local diffeomorphism of $Y$ preserving $Y_0$ followed by  the identification of  pairs $Y_0\subset Y$ with $q(W_h)\subset Z$ induced by the quasi-isomorphism $F$.
\end{remark}

Our techniques establish the following theorem of Strub \cite{Strub}.

\begin{theorem}\label{thm:Strub}
Let $f\colon Y_0\to Z_0$ be a local diffeomorphism  sending $0$ to $0$. Then there is a quasi-isomorphism $F\colon V\to W$.
\end{theorem}
\begin{proof}
The principal points of $Y_0$ are exactly the points where it is locally a smooth manifold without boundary and the codimension one strata are exactly the boundary of the points where $Y_0$ is locally a manifold with boundary. Thus the analogues of Lemmas \ref{lem:principal} and \ref{lem:f_0} hold for $f$ and $f_0$. Since $f_0$ must automatically give an isomorphism of $Y$ and $Z$, the proof of Theorem \ref{thm:main} goes through.
\end{proof}
Finally, we prove a result on lifting global isomorphisms generalizing \cite[Theorem 3.4]{Losik01}.

\begin{theorem}\label{thm;Lifting}
Let $f\colon Y_0\to Z_0 $ be a diffeomorphism. Then there is a lift $F\colon V\to W$ of $f$ where $F$ maps $G$-orbits to $H$-orbits.
\end{theorem}
\begin{proof}
Let $f_t$ be as in Lemma \ref{lem:f_0}. Then  there is a quasi-isomorphism $F'\colon V\to W$ covering $f_0$. Now $f_0\inv f_t$ is an isotopy of the identity of   $Y_0$.  By the isotopy lifting theorem for finite groups \cite{BierstoneLifting} or \cite{SchLifting}, we may lift $f_0\inv f_t$ to an equivariant isotopy $F_t\colon V\times [0,1]\to V$. Then $F' \circ F_1\colon V\to W$ is the desired lift of $f$.
\end{proof}

\section{Lifting biholomorphisms} We have used lifting results from \cite{KrieglLosikMichor03}. The proofs there use results about connections or braid groups. Here we sketch how one can use the slice theorem for finite group actions to obtain their lifting result.

We change notation. Let $G\subset\GL(V)$ and $H\subset\GL(W)$ be finite subgroups where $V$ and $W$ are finite dimensional complex vector spaces.  Let $Y$ denote $V/G$ (the variety associated to $\C[V]^G$) and let $Z$ denote $W/H$.  If $g\in G$ is a pseudoreflection, let $r_g$ denote its order and let $V^{<g>}\subset V$ denote the points of $V^g$ whose isotropy group is generated by $g$. Then $C$, the union of the strata of codimension one, is a disjoint union $C_1\cup\dots\cup C_k$ where each $C_j$ is the image of a subset  $V^{<g_j>}$ for a pseudoreflection $g_j$, $j=1,\dots,k$. To the component $C_j$ we associate the number $r_j:=r_{g_j}$. Similarly we have the union $D$ of the codimension one strata  of $Z$ and a decomposition $D=D_1\cup\dots\cup D_l$ and pseudoreflections $h_1,\dots,h_l\in H$ with orders $s_1,\dots,s_l$ such that $D_j$ is the image of $W^{<h_j>}$, $j=1,\dots,l$.

\begin{theorem} \label {thm:lifting}Suppose that $k=l$ and that 
  $f\colon Y\to Z$ is a biholomorphism which sends each $\overline{C_i}$ onto a $\overline{D_j}$ such that $r_i=s_j$. Then there is a biholomorphic map $F\colon V\to W$ such that $F$ and $F\inv$ send  orbits to orbits, and $F$ induces $f$.
\end{theorem}
\begin{remark}
It is easy to see that the conditions of the theorem are necessary for $F$ to exist.  
\end{remark}
\begin{proof}[Proof of Theorem \ref{thm:lifting}] We use a version of analytic continuation. Let $V_1$ denote the union of the reflection hyperplanes of $V$ and let $V_2$ denote the union of the fixed subspaces $V^K$ of subgroups $K$ of $G$ such that $\codim_V V^K\geq 2$.  Let $Y_j$ denote the image of $V_j$ in $Y$, $j=1$, $2$. Similarly we have subsets $W_j\subset W$ and $Z_j\subset Z$, $j=1$, $2$. Let $v\in V_1\setminus V_2$, let $g$ denote a reflection which generates the isotropy group of $v$ and let $w\in W^h\setminus W_2$ be a point lying above $f(p(v))$ where $h$ is a reflection generating the isotropy group of $w$. By hypothesis, $g$ and $h$ both have the same order, say $r$.  By the slice theorem there is  a ball $B_1$ about $v$ in $V^g$ and a ball $B_2$ in $\C$ around $0$ such that a $G$-neighborhood of $v$ is isomorphic to $G\times^{(g)}(B_1\times B_2)$ where $(g)$ denotes the subgroup generated by $g$. Here $g$ acts trivially on $B_1$ and  acts on $B_2$ via multiplication by a primitive $r$th root of unity. We use the notation $[g_1,x,y]$ to denote the point in $G\times^{(g)}(B_1\times B_2)$ corresponding to $g_1\in G$, $x\in B_1$ and $y\in B_2$. Then the quotient mapping sends a point $[g_1,x,y]$ to $(x,y^r)$ for $x\in B_1$ and $y\in B_2$.  Similarly, there are balls $B_1'$ and $B_2'$ about $w$ in $W^h$ and $0\in \C$ such that the quotient mapping sends $[h_1,x',y')]\in H\times^{(h)}(B_1'\times B_2')$ to $(x',(y')^r)$, $h_1\in H$, $x'\in B_1'$, $y'\in B_2'$.  
The hypotheses on $f$ show that for $(x,y)\in B_1\times B_2^r$, $f(x,y)=(f_1(x,y),f_2(x,y))\in B_1'\times (B_2')^r$ where $f_2$ vanishes when $y=0$ and the derivative of $f_2$ in $y$ has rank 1 along the zero set of $y$. It follows that, locally, $f_2(x,y)$ can be written as $m(x,y)y$ where $m(x,0)$ does not vanish. Thus $f_2(x,y^r)$ has  $r$  holomorphic $r$th  roots along $B_1\times \{0\}$, so that we have $r$  holomorphic lifts   of $f$ in a neighborhood of $v$ which send $v$ to $w$. The lifts are distinguished by their values at any point $(x,y)\in B_1\times B_2$ where $y\neq 0$.

Now let $\gamma(t)$, $0\leq t\leq 1$,  be a continuous curve in $V\setminus V_2$ starting at a base point $v_0$. Let $F_0$ be a germ of a holomorphic map from $V\to W$ which covers $f$. Since $W\setminus W_1\to Z\setminus Z_1$ is a cover, we have an analytic continuation of $F_0$ to $F_t$ at $\gamma(t)$ as long as $\gamma(t)$ never leaves $V\setminus V_1$. However, our argument above shows that we can continue $F$ even if $\gamma(t)$ lands in $V_1\setminus V_2$. Thus we may construct a continuous family $F_t$ of lifts of $f$ along $\gamma(t)$. The family is uniquely determined by $F_0$. Now $V\setminus V_2$ is simply connected, so that $F_1$ only depends upon $\gamma(1)$. Thus we have a lift of $f$ to $V\setminus V_2$. Since $V_2$ has codimension $2$ in $V$, our lift extends to all of $V$. 

We have shown that there is a global lift $F$ of $f$, and similarly there is a global lift $F\inv$ of $f\inv$. Now for any $g\in G$, $F\circ g\circ F\inv$ is an automorphism of $W$ which covers the identity on $W/H$. Thus it must agree with an element of $H$. It follows that $g\mapsto F\circ g\circ F\inv$ gives an isomorphism of $G$ and $H$, and that $F$ and $F\inv$ send orbits to orbits.
\end{proof}


\providecommand{\bysame}{\leavevmode\hbox to3em{\hrulefill}\thinspace}
\providecommand{\MR}{\relax\ifhmode\unskip\space\fi MR }
\providecommand{\MRhref}[2]{%
  \href{http://www.ams.org/mathscinet-getitem?mr=#1}{#2}
}
\providecommand{\href}[2]{#2}

   \end{document}